\documentclass[12pt]{article}

\usepackage{amsthm}
\usepackage{amssymb}
\usepackage{amsmath}
\usepackage{amsfonts}
\usepackage{fullpage}
\usepackage{times}
\usepackage{setspace}

\title{Walk versus Wait:\\The Lazy Mathematician Wins}
\author{\begin{tabular}{ccccc}Justin Gejune Chen\thanks{\texttt{justchen@caltech.edu}} && Scott Duke Kominers\thanks{Corresponding author; \texttt{kominers@fas.harvard.edu}; 8520 Burning Tree Road, Bethesda, MD, 20817.}&& Robert Wyatt Sinnott\thanks{\texttt{rsinnott@fas.harvard.edu}}\\
\small California Institute of Technology&&\small Harvard University&&\small Harvard University\\
\small Pasadena, CA, 91125&&\small Cambridge, MA, 02138&&\small Cambridge, MA, 02138
\end{tabular}}
\date{\today}

\begin{document}
\maketitle

Recently, the first author (Justin) had to walk to the second author (Scott)'s house to work on a problem set.  There is a bus route which covers the distance directly, but the bus arrives sporadically.  So, he faced a conundrum: walk the distance or wait for the bus?  Being lazy,  he would always rather ride the bus, if possible.  Being punctual, however, he will always choose the option which gets him to his destination as quickly as possible.

\paragraph{}Formally, this problem can be stated as follows: Justin has to travel a distance of $d$ miles along a bus route.  (The units are arbitrary, but will ground our discussion.)  Along this route, there are $n$ bus stops $i$, each spaced at a distance of $d_i$ from the starting point.  Without loss of generality, we assume the starting point is the first bus stop, so that $d_1=0$ and that the destination is the last bus stop, so that $d_n=d$.  At each bus stop, Justin is faced with a choice: to walk or to wait.  If he walks on, he can still catch a bus at the next bus stop.

Justin can walk at a velocity of $v_w$ miles per hour and the bus drives at a velocity of $v_b>v_w$ miles per hour.  He must make a decision at the starting point at time $t=0$.  He expects that a bus will arrive at the starting point at time $t$ with probability density function $p(t)$, independently distributed across $t$.

First, we consider the case in which there are only $n=2$ bus stops, the starting point and the destination.  We also temporarily simplify the problem and assume that the next bus arrives at time $t_b>0$, i.e. $p(t_b)=1$ and $p(t)=0$ for $0\leq t<t_b$.  Therefore, the time it would take to walk to the destination is $d/v_w$ hours, while it would take $t_b+d/v_b$ hours to get to the destination on the bus.  The decision is easy, in this case: Justin will walk if and only if \begin{equation}\label{1}\frac{d}{v_w}< t_b+\frac{d}{v_b}.\end{equation}

Even if there are more bus stops between Justin and his destination, his decision in this case is the same.   Indeed, if there is an additional stop at distance $d_2$ from the starting point, then he will never board the bus at the next stop unless, in doing so, he could reach the destination faster than if he walked the distance $d-d_2$.  However, if this is the case, then he would rather save energy and board the bus at the first stop.  It follows that Justin will walk on to the next stop if and only if he stands to gain nothing by waiting for the bus at the first stop, that is, if and only if (\ref{1}) holds.

Weakening the constraint on the probability density function $p(t)$, we find a similar result.  Again, we first consider the case in which  there are only two bus stops, $n=2$.  Justin must minimize his expected travel time, which can be expressed as \begin{equation}\label{2}
\int_0^{t_w}\left(p(t)\left(\frac{d}{v_b}+t\right)\right)dt+\left(1-\int_{0}^{t_w}p(t)\,dt\right)\left(\frac{d}{v_w}+t_w\right),
\end{equation} since the travel time if Justin boards a bus at time $t$ is $\frac{d}{v_b}+t$ and the probability of a bus arriving at time $t$ is $p(t)$.  (We assume that Justin will always take a bus if it arrives.)  Since Justin wants to reach the destination punctually, he sets (\ref{2}) equal to the time required to walk,
\begin{equation}\label{3}\int_0^{t_w}\left(p(t)\left(\frac{d}{v_b}+t\right)\right)dt+\left(1-\int_{0}^{t_w}p(t)\,dt\right)\left(\frac{d}{v_w}+t_w\right)=\frac{d}{v_w}.
\end{equation}

This equation (\ref{3}) implicitly defines an optimal waiting time $t_w$.  If $p(t)$ has an integrable functional form, we can solve for $t_w$ explicitly; for example, if $p(t)=1/t_b$, so that a bus  is believed to arrive by time $t_b$ with a uniform distribution of arrival,  then we find that \begin{align*}
\frac{d}{v_w}&=\int_0^{t_w}\left(\frac{1}{t_b}\left(\frac{d}{v_b}+t\right)\right)dt+\left(1-\int_{0}^{t_w}\frac{1}{t_b}\,dt\right)\left(\frac{d}{v_w}+t_w\right)\\
&=-\frac{t_w^2}{2 t_b}+\left(\frac{d
  \left(\frac{1}{v_b}-\frac{1}{v_w}\right)}{t_b}+1\right)
  t_w+\frac{d}{v_w}.
\end{align*}  We can then solve for $t_w= 2 \left( t_b + \frac{d}{v_b} + \frac{d}{ v_w}\right)$.  Thus in this case waiting  time is dependent on $t_b, d, v_b,$ and $v_w$.

Once again, Justin has nothing to gain by walking, unless he expects to walk the entire distance.  Indeed, if he walks to the second bus stop, then he expects to travel a total of $$\frac{d_2}{v_w} +\int_0^{t_w}\left(\frac{1}{t_b}\left(\frac{d-d_2}{v_b}+t\right)\right)dt+\left(1-\int_{0}^{t_w}p(t)\,dt\right)\left(\frac{d}{v_w}+t_w\right)$$ hours, with $t_w$ implicitly defined by the equation \begin{equation}\int_0^{t_w}\left(\frac{1}{t_b}\left(\frac{d-d_2}{v_b}+t\right)\right)dt+\left(1-\int_{0}^{t_w}p(t)\,dt\right)\left(\frac{d}{v_w}+t_w\right)=\frac{d-d_2}{v_w}.\end{equation}  Substituting, we find that the expected travel time (after a sufficient amount of waiting $t_w$ at the second bus stop) is  $\frac{d}{v_w}$, once again.  Thus, if he intends to wait for the bus at the second stop, Justin could save effort by boarding the bus at the first stop.  It follows inductively that Justin would prefer to walk to the $i$-th stop ($1< i< n$) if and  only if he would prefer to walk all the way to the $n$-th stop.

If Justin must make it to his destination by a fixed time $t_d$, his decision process is nearly unchanged. Letting $$t_{w'}=t_d-\frac{d}{v_w},$$ he will wait at the first bus stop until time $t_w^*:=\min\{t_w,t_{w'}\}$.

\paragraph{}``Eureka!'' Justin shouted upon realizing that the laziest possible waiting strategy would prevail.  He promptly sat down to wait for the bus.

Being a mathematician, Justin had of course chosen to ponder the problem, entirely forgetting the planned meeting.  In the meantime, Scott had grown tired of waiting and had walked over to Justin's house.  When he arrived, he pointed out that it was a holiday, whence no busses were running.

\paragraph{Acknowledgements.} The authors are grateful to Zachary Abel, Cameron Crowe, Todd Kaplan, Paul Kominers, and Stanley Liu  for their helpful comments and suggestions.  They also appreciate Zeeya Merali's assistance in publicizing this article. Additionally, the first two authors are grateful to Mrs.~Susan Schwartz Wildstrom for inspiring and encouraging their interests in mathematics.

\end{document}